\documentclass[11pt]{article}
\usepackage{graphicx, colordvi}
\usepackage{amsfonts}
\usepackage{amssymb}
\usepackage{amsthm,cite,color}
\usepackage{dsfont}
\usepackage{epsfig}
\usepackage{mathrsfs}
\usepackage{amsfonts}
\usepackage{amssymb}
\usepackage{amsmath}
\usepackage{amssymb,amsfonts,amsmath,amsthm,cite,color}
\usepackage{dsfont}
\usepackage{epsfig}
\usepackage{mathrsfs}
\usepackage{longtable}
\usepackage{hyperref}
\hypersetup{
colorlinks=true,
linkcolor=cyan,
filecolor=blue,
urlcolor=red,
citecolor=green,
}

\newtheorem{theo}{Theorem}

\newtheorem{prob}[theo]{Problem}
\newtheorem{coro}[theo]{Corollary}

\newtheorem{lem}[theo]{Lemma}
\newtheorem{exm}[theo]{Example}

\makeatletter \@addtoreset{equation}{section}
\@addtoreset{theo}{section} \makeatother

\newcommand{\bN} { {\mathbb{N}}}

\newcommand{\bQ} { {\mathbb{Q}}}
\newcommand{\bZ} { {\mathbb{Z}}}

\newcommand{\bK} { {\mathbb{K}}}

\newcommand{\lc} { {\mathrm{lc\hspace{0.5ex}}}}

\newcommand{\ann} { {\mathrm{ann\hspace{0.5ex}}}}

\newcommand{\Domb}{\mathrm{Domb}}

\def\qed{\hfill \rule{4pt}{7pt}}
\def\pf{\noindent {\it Proof.} }

\begin{document}
\begin{center}

 {\large \bf Polynomial reduction for holonomic sequences and applications in $\pi$-series and congruences}
\end{center}
\begin{center}
{  Rong-Hua Wang}$^{1}$ and {Michael X.X. Zhong}$^{2}$

   $^1$School of Mathematical Sciences\\
   Tiangong University \\
   Tianjin 300387, P.R. China \\
   wangronghua@tiangong.edu.cn \\[10pt]

   $^2$School of Science\\
   Tianjin University of Technology \\
   Tianjin 300384, P.R. China\\
   zhong.m@tjut.edu.cn
\end{center}

\vskip 6mm \noindent {\bf Abstract.}

Recently, Hou, Mu and Zeilberger introduced a new process of polynomial reduction for hypergeometric terms, which can be
used to prove and generate hypergeometric identities automatically.
In this paper, we extend this polynomial reduction to holonomic sequences.
As applications, we describe an algorithmic way to prove and generate new multi-summation identities.
Especially we present new families of $\pi$-series involving Domb numbers and Franel numbers, and new families of congruence identities for Franel numbers and Delannoy numbers.


\section{Introduction}

For a long time, combinatorial identities, $\pi$-series in particular, were a mysterious part of combinatorics.
It was the seminal work of Wilf and Zeilberger \cite{WZ1990, Zeil1990, Zeilberger1990c, Zeilberger1991} that initiated the study of transforming this mysterious part into science, that everybody, even a computer, could understand.
Since then the mechanical proof of combinatorial identities had received special attention.
Zeilberger's algorithm, also known as the method of creative telescoping, is the core algorithm in the WZ theory.
%Over the past three decades, extensive work has been done around the generalization and application of Zeilberger's algorithm.
%The reduction-based approach is the one which gained much attention as it separates the computation of telescopers and the corresponding certificate and is thus more efficient compared to the original algorithm.

Reduction-based approach plays an indispensable role in the development of the new generation of creative telescoping algorithms, which can separate the calculations of telescopers and certificates for efficiency and construct minimal telescopers.
The first \emph{reduction} algorithm goes back to the work of Ostrogradsky \cite{O1845} and Hermite \cite{Hermite1872} for integrals of rational functions.
In the continuous case, the algorithm was first worked out for bivariate rational functions in \cite{BCCL2010},
and later extended to the multivariate rational case in \cite{BLS2013} using the Griffiths--Dwork method.
The approach has also been extended to algebraic functions \cite{Trager1984,CKK2016}, Fuchsian D-finite functions \cite{CHKK2018}, and general D-finite functions \cite{Hoeven2017,BCLS2018}.

In the discrete case, a reduction-based algorithm was designed for summations of rational functions by Abramov in \cite{Abramov1975}, which was extended to the bivariate rational case in \cite{CHHLW}.
The hypergeometric and holonomic cases were studied by Abramov and Petkov{\v{s}}ek in \cite{AP2001,AP2002} and van der Hoeven in \cite{Hoeven2018} respectively. In 2015, Chen et al. \cite{CHKL2015} introduced the concept of \emph{polynomial reduction} in the modified Abramov--Petkovcek algorithm which is more efficient and can be used to compute minimal telescopers for bivariate hypergeometric terms.

There are two fundamental problems related to telescopers:
one is deciding whether telescopers exist for a given function, the other one is designing efficient algorithms for constructing telescopers when exist.
%The Ostrogradsky--Hermite reduction plays a crucial role in solving both problems.
For the existence problem, bivariate mixed hypergeometric case can be solved via the Ostrogradsky--Hermite reduction \cite{CCFFL2015}, while
trivariate mixed rational case was handled by the extended reduction in \cite{CHLW2016,CDZ2019}.
For the construction problem, the Ostrogradsky--Hermite reduction and its variants have been applied in \cite{CHKL2015,Huang2016,CKK2016,CHKK2018,Hoeven2017,Hoeven2018} for bivariate functions,
and in \cite{BLS2013,Lairez2016,CHHLW} beyond the bivariate case.

Although the method of creative telescoping is a powerful tool in the mechanical proof of combinatorial identities, the reduction itself was rarely used directly in the verifications.
In 2021, Hou, Mu and Zeilberger~\cite{HouMuZeil2021} introduced another \emph{polynomial reduction} process, which can be employed to derive infinite families of supercongruences \cite{HouMuZeil2021} and new hypergeometric identities \cite{HouLi2021}.
Recently, Wang and Zhong \cite{WZ2022} generalized the polynomial reduction to the $q$-rational case.
This makes it possible to prove and discover $q$-identities automatically.
Especially several $q$-analogues of Ramanujan-type series for $\pi$ were presented.

The Hou--Mu--Zeilberger reduction and its variants were all designed for ($q$-)hypergeometric terms.
In this paper, we generalize the Hou--Mu--Zeilberger reduction to the holonomic case.
This enables us to handle multi-summation identities.
As applications, we provide an algorithmic way to prove and discover new series for $\pi$ involving Domb numbers and Franel numbers.
This will confirm and generalize some of Z.-W Sun's conjectures in \cite{Sun2021}.
New families of congruence identities on Franel numbers and Delannoy numbers are also obtained.

\section{Polynomial reduction for holonomic sequences}

Let $\bK$ be a field of characteristic $0$.
A sequence $(F(n))_{n=0}^{\infty}$ is \emph{holonomic} over $\bK$ if there exist polynomials $a_0(n),a_1(n),\ldots,a_J(n)\in\bK[n]$ with $a_J(n)\neq 0$ such that
\begin{equation}\label{eq:rec homo}
\sum_{i=0}^{J}a_i(n)F(n+i)=0.
\end{equation}
Or, equivalently, if we define the \emph{annihilator} of $F(n)$ by
\[
\ann F(n):=\left\{L=\sum_{i=0}^{J}a_i(n)\sigma^i\in \bK[n][\sigma]\mid L(F(n))=0\right\},
\]
where $\sigma$ is the shift operator (that is, $\sigma F(n)=F(n+1)$),
then $(F(n))_{n=0}^{\infty}$ is holonomic if and only if $\ann F(n)\neq \{0\}$.
We call $J$ in \eqref{eq:rec homo} the order of the recurrence relation for $F(n)$, and the minimum order of all such recurrences is called the \emph{order} of $F(n)$.

The class of holonomic sequences covers a great percentage of combinatorial sequences arising in applications.
For example, harmonic numbers, Fibonacci numbers, Domb numbers, Franel numbers and all hypergeometric sequences are holonomic.
Specifically, given a hypergeometric term $t_n$, there exist polynomials $a(n),b(n)\in\bK[n]$ such that
\[\frac{t_{n+1}}{t_n}=\frac{a(n)}{b(n)},\]
that is,
\[a(n)t_n-b(n)t_{n+1}=0.\]
In the polynomial reduction introduced by Hou, Mu and Zeilberger \cite{HouMuZeil2021}, a key step is to characterize such polynomials $p(n)\in\bK[n]$ that the product $p(n)t_n$ is \emph{Gosper-summable}, that is,
\[p(n)t_n=\Delta (u(n)t_n),\]
for some rational function $u(n)\in\bK(n)$, where $\Delta$ is the difference operator (that is, $\Delta F(n)=\sigma F(n)-F(n)=F(n+1)-F(n)$).

It is natural to consider a similar problem in the holonomic case.
\begin{prob}
Given a holonomic sequence $(F(n))_{n=0}^{\infty}$ satisfying \eqref{eq:rec homo}, for which polynomials $q(n)\in\bK[n]$, the product $q(n)F(n)$ can be written as
\begin{equation}\label{eq: Abramov--van-Hoeij}
q(n)F(n)=\Delta\left(\sum_{i=0}^{J-1}u_i(n)F(n+i)\right)
\end{equation}
for some rational functions $u_0(n),u_1(n),\ldots,u_{J-1}(n)\in\bK(n)$?
\end{prob}

For any operator $L=\sum_{i=0}^{J}a_i(n)\sigma^i$ with $a_i(n)\in\bK[n]$,
the \emph{adjoint} of $L$ is defined by
\begin{equation}\label{eq:L}
L^{\ast}=\sum_{i=0}^{J}\sigma^{-i}a_i(n).
\end{equation}
Then for any polynomial $p(n)\in\bK[n]$,
\begin{equation*}\label{eq:B}
L^{\ast}(p(n))=\sum_{i=0}^{J}a_i(n-i)p(n-i).
\end{equation*}

In 2018, van der Hoeven \cite[Proposition 3.2]{Hoeven2018} derived the following difference Lagrange identity
\begin{equation}\label{eq:Lagrange identity}
p(n)L(F(n))-L^{\ast}(p(n))F(n)=\Delta
\left(\sum_{i=0}^{J-1}u_i(n)F(n+i)\right)
\end{equation}
by induction,
where
\begin{equation}\label{eq:u}
u_i(n)=\sum_{j=1}^{J-i}a_{i+j}(n-j)p(n-j).
\end{equation}
Next we provide another proof of \eqref{eq:Lagrange identity} by a direct calculation and the fact
\begin{equation}\label{eq:summable reduction}
(\sigma^i-1)g(n)=\Delta\left(\sum_{j=0}^{i-1}g(n+j)\right),\quad \forall i>0.
\end{equation}

\begin{lem}\label{lem: summable}
Let $L=\sum_{i=0}^{J}a_i(n)\sigma^i$ with $a_i(n)\in\bK[n]$.
Then \eqref{eq:Lagrange identity} holds for any polynomial $p(n)\in\bK[n]$.
\end{lem}
\pf
By the definition of $L^{\ast}$ and equality \eqref{eq:summable reduction}, we have
\begin{align*}
 &p(n)L(F(n))-L^{\ast}(p(n))F(n)\\
%=&\sum_{i=1}^{J}\left(a_i(n)p(n)F(n+i)-a_i(n-i)p(n-i)F(n)\right)\\
=&\sum_{i=1}^{J}(\sigma^i-1)(a_i(n-i)p(n-i)F(n))\\
=&\Delta\left(\sum_{i=1}^{J}\sum_{j=0}^{i-1}a_i(n-i+j)p(n-i+j)F(n+j)\right)\\
=&\Delta\left(\sum_{j=0}^{J-1}\sum_{i=j+1}^{J}a_i(n-i+j)p(n-i+j)F(n+j)\right)\\
=&\Delta\left(\sum_{i=0}^{J-1}\sum_{j=i+1}^{J}a_j(n-j+i)p(n-j+i)F(n+i)\right)\\
=&\Delta\left(\sum_{i=0}^{J-1}u_i(n)F(n+i)\right).
\end{align*}
\qed

When $L\in\ann F(n)$, identity \eqref{eq:Lagrange identity} reduces to
\begin{equation}\label{eq: summable}
L^{\ast}(p(n))F(n)=\Delta
\left(-\sum_{i=0}^{J-1}u_i(n)F(n+i)\right).
\end{equation}
Hence $q(n)=L^{\ast}(p(n))$
is a desired polynomial in $\bK[n]$ such that \eqref{eq: Abramov--van-Hoeij} holds.
From equality \eqref{eq: summable} we obtain
\begin{equation}\label{eq:rec ann}
\sum_{k=0}^{n-1}L^{\ast}(p(k)) F(k)
=\left(\sum_{i=0}^{J-1}u_i(0)F(i)\right)
-\left(\sum_{i=0}^{J-1}u_i(n)F(n+i)\right),
\end{equation}
where $u_i$ is defined in \eqref{eq:u}.
Taking $n\to\infty$, we usually get
\[\sum_{n=0}^{\infty}L^{\ast}(p(n))F(n)
=C,
\]
where $C$ is a constant.

From the proof of Lemma \ref{lem: summable}, one can see  \eqref{eq:Lagrange identity} and \eqref{eq: summable} still hold for any rational function $p(n)\in\bK(n)$.
When the order $J$ in \eqref{eq:rec homo} is minimum, the Abramov--van-Hoeij algorithm \cite{AbramovvanHoeij,KauersPaule2011}
ensures that all rational functions $q(n)\in\bK(n)$ such that \eqref{eq: Abramov--van-Hoeij} holds are of the form $L^{\ast}(p(n))$ with $p(n)\in\bK(n)$.

In this paper, to make the reduction work, $L^{\ast}(p(n))$ needs to be a polynomial.
Abramov \cite{Abramov1989,Abramov1995} characterized the denominator of a rational solution $p(n)$ to a difference equation in the form
\[
a_0(n)p(n)+a_1(n-1)p(n-1)+\cdots+a_J(n-J)p(n-J)=b(n),
\]
where $a_i(n)\in\bK[n]$, $0\leq i\leq J$, and $b(n)\in\bK[n]$.
It may happen that $p(n)\in\bK(n)\setminus\bK[n]$ but $L^{\ast}(p(n))\in\bK[n]$.
The following lemma shows that this rarely happens.

\begin{lem}\label{lem:rational->polynomial}
Let $L=\sum_{i=0}^{J}a_i(n)\sigma^i$ with $a_i(n)\in \bK[n],\ 0\leq i\leq J$ and $a_0(n)a_J(n)\neq 0$.
If $p(n)\in\bK(n)$ and
\begin{equation}\label{eq:coprime}
\gcd (a_0(n),a_J(n+i))=1,\quad\forall i\in\bN,
\end{equation}
then
$L^{\ast}(p(n))$
is a polynomial in $\bK[n]$ if and only if $p(n)\in\bK[n]$.
\end{lem}

\pf
Suppose $p(n)=\frac{r(n)}{s(n)}$ with $r(n),s(n)\in\bK[n]$ and $\gcd(r(n),s(n))=1$.
If $q(n)=L^{\ast}(p(n))$ is a polynomial in $\bK[n]$, then
\begin{align*}
q(n)
=\sum_{i=0}^{J}a_i(n-i)p(n-i)
=\sum_{i=0}^{J}a_i(n-i)\frac{r(n-i)}{s(n-i)}.
\end{align*}
Multiplying $s(n)s(n-1)\cdots s(n-J)$ on both ends, we obtain
\begin{align*}
q(n)\prod_{i=0}^{J}s(n-i)
=\sum_{i=0}^{J}a_i(n-i)r(n-i)\prod_{\substack{0\leq j\leq J\\j\neq i}} s(n-j).
\end{align*}
Apparently, $s(n)$ is a divisor of the left hand side, then we must have
\begin{align*}
s(n)\mid a_0(n)r(n)s(n-1)\cdots s(n-J).
\end{align*}
From $\gcd(r(n),s(n))=1$ we know that
\begin{equation}\label{eq: division1}
s(n)\mid a_0(n)s(n-1)\cdots s(n-J).
\end{equation}
By a similar argument,
\[s(n-J)\mid a_J(n-J)s(n)s(n-1)\cdots s(n-J+1),\]
namely,
\begin{equation}\label{eq: division2}
s(n)\mid a_J(n)s(n+J)s(n+J-1)\cdots s(n+1).
\end{equation}

Suppose that $s(n)$ is not a constant and $t(n)$ is an irreducible factor of $s(n)$.
Then by \eqref{eq: division1}, if $t(n)\nmid a_0(n)$, then $t(n)\mid s(n-j_1)$ for some $j_1>0$, that is, $t(n+j_1)\mid s(n)$.
Again, by \eqref{eq: division1}, if $t(n+j_1)\nmid a_0(n)$, then $t(n+j_1+j_2)\mid s(n)$ for some $j_2>0$.
Since $s(n)$ can not have infinitely many distinct factors, there must exist a $j\ge 0$ such that $t(n+j)\mid a_0(n)$.
By a similar argument, \eqref{eq: division2} guarantees that there exists an $i\ge 0$ such that $t(n-i)\mid a_J(n)$.
Then $\gcd(a_0(n),a_J(n+i+j))$ is not a constant, contradicting \eqref{eq:coprime}.
So $s(n)$ must be a constant, namely, $p(n)$ is a polynomial.
As the converse is clearly true, this completes the proof.
\qed

Combining Lemma \ref{lem: summable} and Lemma \ref{lem:rational->polynomial}, we have the following result.

\begin{theo}\label{thm:rational to polynomial}
Let $L=\sum\limits_{i=0}^{J}a_i(n)\sigma^i\in\ann F(n)$ with  $a_0(n)a_J(n)\neq 0$ and $J>0$ the order of $F(n)$.
If
\begin{equation*}
\gcd (a_0(n),a_J(n+i))=1,\quad\forall i\in\bN,
\end{equation*}
then $q(n)\in\bK[n]$ satisfies
\eqref{eq: Abramov--van-Hoeij} if and only if
$q(n)=L^{\ast}(p(n))$ for some $p(n)\in\bK[n]$.
\end{theo}

Now we assume that $p(n)\in\bK[n]$ and try to determine the degree of $L^{\ast}(p(n))$ when $L=\sum\limits_{i=0}^{J}a_i(n)\sigma^i$ is given.
To this aim, some notations are needed.
Let
\begin{equation}\label{eq:d and bk}
b_k(n)=\sum_{j=k}^{J}\binom{j}{k}a_{J-j}(n+j-J)\text{ and }
d=\max_{0\leq k \leq J} \{\deg b_k(n)-k\}.
\end{equation}
Note that
\[
f(s)=\sum_{k=0}^{J}[n^{d+k}](b_k(n))s^{\underline{k}}
\]
is a nonzero polynomial in $s$.
Here $[n^{d+k}](b_k(n))$ denotes the coefficient of $n^{d+k}$ in $b_k(n)$ and $s^{\underline{k}}$ denotes the falling factorial defined by $s^{\underline{k}}=s(s-1)\cdots (s-k+1)$.
Let
\begin{equation}\label{eq:nonnegative roots}
R_{L}=\{s\in\bN \mid f(s)=0\}.
\end{equation}
Then $L$ is called \emph{degenerated} if $R_{L}\neq \emptyset$.

\begin{lem}\label{lem:degree}
Let $L=\sum\limits_{i=0}^{J}a_i(n)\sigma^i$ and $d$ be given by \eqref{eq:d and bk}.
Then for any nonzero polynomial $p(n)$, we have
\[
\deg L^{\ast}(p(n))\left\{
     \begin{array}{ll}
       < d+\deg p(n), & \hbox{if $L$ is degenerated and $\deg p(n)\in R_{L}$,} \\
       =d+\deg p(n), & \hbox{otherwise.}
     \end{array}
   \right.
\]
\end{lem}
\proof
Let $q(n)=p(n-J)$.
Notice that $\sigma=E+\Delta$, where $E$ is the identity map.
Then
\begin{align}\label{eq: L(p(n))}
  L^{\ast}(p(n))
& = \sum_{i=0}^{J}a_i(n-i)q(n+J-i)
 = \sum_{i=0}^{J}a_i(n-i)(E+\Delta)^{J-i}q(n)\nonumber \\
& = \sum_{j=0}^{J}a_{J-j}(n+j-J)(E+\Delta)^{j}q(n)\nonumber \\
&= \sum_{j=0}^{J}a_{J-j}(n+j-J)\sum_{k=0}^{j}\binom{j}{k}\Delta^k (q(n))\nonumber \\
&=\sum_{k=0}^{J}b_k(n)\Delta^k(q(n)),
\end{align}
where $b_k(n)$ is defined in \eqref{eq:d and bk}.

Denote $s:=\deg q(n)=\deg p(n)$.
Then $\Delta^k(q(n))$ is a polynomial of degree $s-k$ and $\lc(\Delta^k(q(n)))=\lc(q(n)) s^{\underline{k}}$ if $k\le s$,
and $\Delta^k(q(n))=0$ if $k>s$.
Here $\lc(f(n))$ denotes the leading coefficient of the polynomial $f(n)$.
By equality \eqref{eq: L(p(n))}, we know
\[
\deg L^{\ast}(p(n)) \leq \max_{0\leq k\leq J}\{\deg b_k(n)+\deg q(n)-k\}=d+s.
\]
Noting that $s^{\underline{k}}=0$ if $k>s$, it is easy to see that $\deg L^{\ast}(p(n)) <d+s$ if and only if
$
\sum_{k=0}^{J}[n^{d+k}](b_k(n))s^{\underline{k}}=0,
$
which means $L$ is degenerated and $s=\deg p(n)\in R_{L}$.
\qed

With this lemma, we are able to give a precise description of the polynomial reduction process for holonomic sequences.

\emph{The polynomial reduction process}:
Let $L=\sum\limits_{i=0}^{J}a_i(n)\sigma^i\in\bK[n][\sigma]$ with $a_J\neq 0$ and
\begin{equation}\label{eq:p_s}
  q_s(n)=L^{\ast}(p_s(n))=\sum_{i=0}^{J}a_i(n-i)p_s(n-i),
\end{equation}
where $p_s(n)$ is a polynomial in $\bK[n]$ of degree $s\in\bN$.
We first consider the case when $L$ is not degenerated.
By Lemma \ref{lem:degree}, we know
\[
\deg q_s(n)=d+s,\quad \forall s\in\bN,
\]
where $d$ is defined as \eqref{eq:d and bk}.
Then for any polynomial $Q(n)$ of degree $m$ with $m\geq d$,
it can be written by the division algorithm as
\begin{equation}\label{eq:reduction step}
 Q(n)=\sum_{s=0}^{m-d}c_{s}q_{s}(n)+\tilde{q}(n),
\end{equation}
where $c_s\in\bK$ for $0\leq s\leq m-d$ and $\tilde{q}(n)$ is a polynomial of degree less than $d$.
When $L$ is degenerated, by Lemma \ref{lem:degree},
\[
\deg q_s(n)=d+s,\quad \forall s\in \bN\setminus R_{L}.
\]
Then \eqref{eq:reduction step} works well except for the polynomials of degree $d+s$ for $s\in R_{L}$.
Thus for any polynomial $Q(n)$ of degree $m$ with $m\geq d$, we can write it as
\begin{equation}\label{eq:reduction step1}
 Q(n)=
 \sum_{\substack{0\leq s\leq m-d\\s\notin R_{L}}}c_{s}q_{s}(n)
 +\sum_{\substack{0\leq s\leq m-d\\s\in R_{L}}}c_{s}n^{d+s}+\tilde{q}(n),
\end{equation}
where $c_s\in\bK$ for $0\leq s\leq m-d$ and $\tilde{q}(n)$ is a polynomial with $\deg\tilde{q}(n)<d$.
Equality \eqref{eq:reduction step} (or \eqref{eq:reduction step1}) is called the \emph{polynomial reduction} with respect to $L$ when it is not degenerated (or degenerated).

%In the next section, we will show how to use the polynomial reduction to derive new identities with holonomic terms from the old.
\section{Applications}
In this section,  we will take Domb numbers, Franel numbers and Delannoy numbers as examples to illustrate how to generate new $\pi$-series and congruence identities algorithmically by the polynomial reduction.
\subsection{\texorpdfstring{Generating new $\pi$-series}{Generating new π-series}}
The \emph{Domb numbers} are given by
\[
\Domb(n)=\sum_{k=0}^n \binom{n}{k}^2\binom{2k}{k}\binom{2(n-k)}{n-k}.
\]
Chan, Chan and Liu \cite{CCL2004} and Rogers \cite{Rogers2009} derived
\begin{equation}\label{eq:Domb -32}
\sum_{n=0}^{\infty}(5n+1)\frac{\Domb(n)}{64^n}=
\frac{8\sqrt{3}}{3 \pi}
\text{ and }
\sum_{n=0}^{\infty}(3n+1)\frac{\Domb(n)}{(-32)^n}=
\frac{2}{\pi}.
\end{equation}

The following identity was conjectured by Z.-W. Sun \cite{Sun2021}.
We will take it as an example to show how to use the polynomial reduction method to prove new identities from the old.

\begin{theo}
\begin{equation}\label{eq:-32 and k^2}
\sum_{n=0}^{\infty}n^2(n-1)(9n+1)\frac{\Domb(n)}{(-32)^n}=
\frac{4}{3\pi}.
\end{equation}
\end{theo}

\pf Let $F(n)=\Domb(n)/(-32)^n$.
By Zeilberger's algorithm, we find that
$
L=\sum_{i=0}^{2} a_i(n) \sigma^i\in\ann F(n),
$
with
$a_0(n)=(n+1)^3,a_1(n)=(2n+3) (5 n^2+15n+12)$ and
$a_2(n)= 16 (n+2)^3$.
Then it is easy to see that
\[d=3\text{ and }R_L=\emptyset.\]
So $L$ is nondegenerated.
By Lemma \ref{lem:degree}, $\deg L^{\ast}(p_s(n))=s+3$ for any polynomial $p_s(n)\in\bK[n]$ of degree $s\ge 0$.
Substituting $a_0(n),a_1(n),a_2(n)$ and $F(n)$ into equality \eqref{eq:rec ann} and taking $n\to \infty$ leads to
\begin{equation*}\label{eq:final equ}
\sum_{n=0}^{\infty} L^{\ast}(p_s(n))F(n)=0.
\end{equation*}
for any polynomial $p_s(n)$.

The polynomial reduction shows that
\[
n^2(n-1)(9n+1)=\frac{2}{3}(3n+1)+\frac{1}{3}L^{\ast}(n).
\]
Multiplying by $\frac{\Domb(n)}{(-32)^n}$ on both sides of the above identity and then summing over $n$ from $0$ to $\infty$, we derive
\[
\sum_{n=0}^{\infty}n^2(n-1)(9n+1)\frac{\Domb(n)}{(-32)^n}=
\frac{2}{3}\sum_{n=0}^{\infty}(3n+1)\frac{\Domb(n)}{(-32)^n}
=\frac{4}{3\pi}
\]
with the help of identity \eqref{eq:Domb -32}.
This completes the proof of \eqref{eq:-32 and k^2}.
\qed

The above theorem shows how to algorithmically prove a conjectured identity from a known one.
The method can actually confirm all those identities listed in Conjecture 8.2 of \cite{Sun2021} by Z.-W. Sun.
However, to pose these conjectures needs genuine intuition, insight, experience and hard work.

For the rest part of this paper, we are going to give an algorithmic way of generating new identities from the old, which does not depend on intuition or experience at all.
For illustration, we take the following $\pi$-series involving Domb numbers as a starting point:
\begin{equation}\label{eq:Domb}
\sum_{n=0}^{\infty}(un+v)\frac{\Domb(n)}{m^n}=
\frac{\lambda}{\pi}.
\end{equation}
\begin{theo}\label{th: Domb general identities}
Suppose an identity of the form \eqref{eq:Domb} holds for some nonzero $m\in\bZ$ and $u,v\in\bQ$.
Then for any nonconstant polynomial $P(n)\in\bQ[n]$, one can find a nonzero polynomial $Q(n)\in\bQ[n]$ with $\deg Q(n)\leq 2$ and a constant $c\in\bQ$ such that
\begin{equation}\label{eq:Domb -32 generalized}
\sum_{n=0}^{\infty}P(n) Q(n)\frac{\Domb(n)}{m^n}=\frac{c\lambda}{\pi}.
\end{equation}
\end{theo}
\pf Let $F(n)=\Domb(n)/m^n$.
By Zeilberger's algorithm, we find that
$
L=\sum_{i=0}^{2} a_i(n) \sigma^i\in\ann F(n),
$
with
$a_0(n)=64(n+1)^3,a_1(n)=-2m(2n+3)(5 n^2+15n+12)$ and
$a_2=m^2(n+2)^3$.
Let $d, R_L$ be defined by \eqref{eq:d and bk} and \eqref{eq:nonnegative roots}.
Then it is easy to see that ,
\[d=3\text{ and }R_L=\emptyset  \text{ for } m\notin \{4,16\}\]
while
\[d=2\text{ and }R_L=\emptyset  \text{ for } m\in \{4,16\}.\]
So $L$ is nondegenerated.
Then by Lemma \ref{lem:degree}, $\deg L^{\ast}(p_s(n))=s+d$ for any polynomial $p_s(n)\in\bK[n]$ of degree $s\ge 0$.

%Note that $a_2(-2)=0$ and $a_1(-1)F(0)+a_2(-1)F(1)=0$.
Substituting the above $a_0(n),a_1(n),a_2(n)$ and $F(n)$ into equality \eqref{eq:rec ann} and taking $n\to \infty$ leads to
\begin{equation*}
\sum_{n=0}^{\infty} L^{\ast}(p_s(n))F(n)
=0
\end{equation*}
for any polynomial $p_s(n)$ of degree $s$.
In the following, we take $p_s(n)=n^s$.

For any non-constant polynomial $P(n)$, let $\ell=\deg P(n)>0$.
Suppose
$
Q(n)=e_0+e_1n+e_2n^2
$
with indeterminants $e_i,\ i=0,1,2$.
Now solve the equation
\begin{equation}\label{eq:general equ}
P(n)Q(n)=c(un+v)+c_0 L^{\ast}(n^0)+c_1 L^{\ast}(n^1)\cdots +c_{\ell} L^{\ast}(n^{\ell})
\end{equation}
for indeterminats $e_0,e_1,e_2,c,c_0,\ldots,c_{\ell}$ in $\bQ$.
By comparing the coefficients of $n^k$ on both sides for $k=0,1,\ldots,\ell+3$, we get a system of $\ell+4$ homogeneous linear equations in $\ell+5$ indeterminants,
so there must be nonzero solutions for $e_0,e_1,e_2,c,c_0,\ldots,c_{\ell}$.
Notice that $\deg L^{\ast}(n^s)=s+d$.
Therefore, $e_0,e_1,e_2$ can not be all zero, that is, there is a nonzero polynomial $Q(n)\in\bQ[n]$ with $\deg Q(n)\le 2$ such that \eqref{eq:general equ} holds.

Multiplying by $F(n)$ on both sides of \eqref{eq:general equ} and then summing over $n$ from $0$ to $\infty$, we obtain
\[
 \sum_{n=0}^{\infty} P(n)Q(n)\frac{\Domb(n)}{m^n}
 = c\sum_{n=0}^{\infty} (un+v)\frac{\Domb(n)}{m^n}
 =\frac{c\lambda}{\pi}.
\]
%This completes the proof.
\qed

Using the method in the proof of Theorem \ref{th: Domb general identities}, one can not only confirm Z.-W. Sun's Conjecture 8.2 in \cite{Sun2021}, but also generate as many new ones as you like.

As an example, we now revive the discovery of \eqref{eq:-32 and k^2} with our method.
Take $m=-32$, one can see
$L=\sum_{i=0}^{2} a_i(n) \sigma^i\in\ann \frac{\Domb(n)}{(-32)^n}$, where
$a_0(n)=(n+1)^3,a_1(n)=(2n+3) (5 n^2+15n+12)$ and
$a_2(n)= 16 (n+2)^3$.
Let $P(n)=n^2$ and solve
\begin{equation}\label{eq:-32 ex}
P(n)(e_0+e_1n+e_2n^2)=c(3n+1)+c_0 L^{\ast}(n^0)+c_1 L^{\ast}(n^1).
\end{equation}
We find for any $c\in\bQ,$
\[
c_0=0, c_1=c/2,
e_0= -3 c/2, e_1= -12 c,
e_2=27 c/2
\]
is a solution of \eqref{eq:-32 ex}.
Taking $c=2/3$, we arrive at
\[
n^2(n-1)(9n+1)=\frac{2}{3}(3n+1)+\frac{1}{3}L^{\ast}(n^1).
\]
Multiplying by $\frac{\Domb(n)}{(-32)^n}$ on both sides of the above identity and then summing over $n$ from $0$ to $\infty$, we derive
\eqref{eq:-32 and k^2}.

Different choices of $P(n)$ may lead to different identities,
for example, we obtain
\[
\sum_{n=0}^{\infty}(n^2+n+1)(126n^2+41n+5)\frac{\Domb(n)}{(-32)^n}=
-\frac{100}{3\pi}
\]
by taking $P(n)=n^2+n+1$.

From the discussion above, one can see how the polynomial reduction may be applied to holonomic sequences.
Here are more examples.

The \emph{Franel numbers} and \emph{Franel numbers of order $4$} are defined respectively by
\begin{equation*}\label{eq:franel}
f_n=\sum_{k=0}^{n} \binom{n}{k}^3 \text{and }
f_n^{(4)}=\sum_{k=0}^{n} \binom{n}{k}^4.
\end{equation*}
Many series for $\pi$ involving $f_n$ and $f_n^{(4)}$ are obtained via modular forms in \cite{Cooper2012, CTYZ2011}.
Those series are of the form
\begin{equation}\label{eq:Franel}
 \sum_{n=0}^{\infty} (un+v)\frac{A(n)}{m^n}=\frac{\lambda\sqrt{\alpha}}{\pi},
\end{equation}
where $A(n)=f_n^{(4)}$ or $A(n)=\binom{2n}{n}f_n$, $u,v,m\in\bZ$, $\lambda\in\bQ$ with $\lambda um\neq 0$ and $\alpha$ is a positive integer.

One can check that $A(n)$ satisfies a recurrence relation of order $2$.
Similar to the proof of Theorem \ref{th: Domb general identities}, we obtain the following result.
\begin{theo}\label{th: franel general identities}
Suppose an identity of the form \eqref{eq:Franel} holds for some $u,v,m\in\bZ$ with $m\neq 0$.
Then for any nonconstant polynomial $P(n)\in\bQ[n]$, one can find a nonzero polynomial $Q(n)\in\bQ[n]$ with $\deg Q(n)\leq 2$ and a constant $c\in\bQ$ such that
\begin{equation}\label{eq:franel generalized}
\sum_{n=0}^{\infty}P(n) Q(n)\frac{\text{A}(n)}{m^n}=\frac{c\lambda\sqrt{\alpha}}{\pi}.
\end{equation}
\end{theo}
%The proof of Theorem \ref{th: franel general identities} can be illustrated through the following example.
\begin{exm}
We can prove the following conjectural identity by Z.-W. Sun\cite[Conjecture 8.3]{Sun2021}
\begin{align}\label{eq:fk5776}
\sum_{n=1}^{\infty}
&n^3(47808294003072n^2-102482715691400n
+52422407372915)\frac{f_n^{(4)}}{5776^n}\nonumber\\
&=-\frac{122626206796\sqrt{95}}{625\pi}
\end{align}
utilizing the polynomial reduction and the identity
\[
\sum_{n=1}^{\infty}
(408n+47)\frac{f_n^{(4)}}{5776^n}
=\frac{1444\sqrt{95}}{95\pi}
\]
derived by Cooper~\cite{Cooper2012}.

Taking $m=5776$, one can see
$L=\sum_{i=0}^{2} a_i(n) \sigma^i\in\ann\frac{f_n^{(4)}}{5776^n}$, where
$a_0(n)=-(n+1)(4n+3)(4n+5),
a_1(n)=-2888(2n+3)(3n^2+9n+7)$ and
$a_2(n)= 8340544 (n+2)^3$.
Let $P(n)=n^3$ and solve
\begin{equation}\label{eq:5776 ex}
P(n)(e_0+e_1n+e_2n^2)=c(3n+1)+c_0 L^{\ast}(n^0)+c_1 L^{\ast}(n^1)+c_2 L^{\ast}(n^2).
\end{equation}
For any $c\in\bQ$, let $d=c/1613502721$.
One can check that
$c_0=-590794567 d,$
$ c_1=-1338119121 d,$
$c_2=-717997495 d,$
$e_0= -6552800921614375 d,$
$e_1= 12810339461425000 d,$
$e_2=-5976036750384000 d$
is a solution of \eqref{eq:5776 ex}.
Take $c=-1613502721/125$.
Multiplying by $\frac{f_n^{(4)}}{5776^n}$ on both sides of \eqref{eq:5776 ex} and then summing over $n$ from $0$ to $\infty$, we derive \eqref{eq:fk5776}.
\end{exm}

Similarly, we can confirm all those identities in Conjecture 8.3 and Conjecture 8.4 (\romannumeral1) of Z.-W. Sun's paper \cite{Sun2021}.

\subsection{Generating new congruence identities}
In this subsection, we will show the polynomial reduction method can also be applied to prove and discover new families of congruence identities.
We will take the Franel numbers $f_k$ and the central Delannoy numbers $D_k$ as examples to reveal the process.

In 2013, Z.-W. Sun \cite{Sun2013} initiated the systematic investigation of fundamental congruences of Franel numbers.
Many interesting congruences are obtained, for example, for any prime $p>3$ there hold
\begin{align}
  &\sum_{k=0}^{p-1}(-1)^k f_k \equiv \left(\frac{p}{3}\right)\pmod {p^2}, \label{eq:franel 0}\\
  &\sum_{k=0}^{p-1}k(-1)^k f_k \equiv -\frac{2}{3}\left(\frac{p}{3}\right)\pmod {p^2} , \label{eq:franel 1}\\
  &\sum_{k=0}^{p-1}k^2(-1)^k f_k \equiv \frac{10}{27}\left(\frac{p}{3}\right)\pmod {p^2}. \label{eq:franel 2}
\end{align}
Here $\left(\dfrac{a}{p}\right)$ denotes the Legendre symbol.
Later V.J.W. Guo \cite{Guo2013} confirmed the following two conjectures by Z.-W. Sun \cite{Sun2013},
\begin{align}
& \sum_{k=0}^{n-1}(3k+2)(-1)^k f_k\equiv 0 \pmod{2n^2},\label{eq:franel 3k+2}\\
& \sum_{k=0}^{p-1}(3k+2)(-1)^k f_k\equiv 2p^2(2^p-1)^2 \pmod{p^5}.\label{eq:franel 3k+2 mod p}
\end{align}
Recently, Wang and Sun \cite{WangSun2019} derived more divisibility results on Franel numbers like
\begin{align}
  &9\sum_{k=1}^{n}k^2(3k+1)(-1)^k f_k \equiv 0\pmod {2n^2(n+1)^2}, \label{eq:franel 3-1}\\
  &3\sum_{k=1}^{n}(9k^3-15k^2-10k)(-1)^k f_k \equiv 0\pmod {4n(n+1)^2} \label{eq:franel 3-2}.
\end{align}
In 2021, by telescopings of $P$-recursive sequences, Hou and Liu \cite{HouLiu2021} found
\begin{equation}\label{eq:franel 3k+2 D}
3\sum_{k=0}^{n-1}(3k+2)(-1)^k f_k=
n^2((-1)^n f_n+8(-1)^{n-1} f_{n-1}),
\end{equation}
which reproves \eqref{eq:franel 3k+2} since
\begin{align}
& (-1)^nf_n\equiv \sum_{k=0}^{n}\binom{n}{k}
\equiv 2^n \equiv 0\pmod 2, n\geq 1, \label{eq:franel mod2}\\
& (-1)^nf_n\equiv (-1)^n\sum_{k=0}^{n}\binom{n}{k}
\equiv (-2)^n \equiv 1\pmod 3.\nonumber
\end{align}
When $n=p>3$ is a prime in \eqref{eq:franel 3k+2 D}, direct calculations lead to equality \eqref{eq:franel 3k+2 mod p} since $f_p\equiv 2 \pmod {p^3}$ and
\[
f_{p-1}\equiv 1+3(2^{p-1}-1)+3(2^{p-1}-1)^2 \pmod {p^3},
\]
as proved by Z.-W. Sun \cite{Sun2013}.

Next, we will generalize \eqref{eq:franel 3k+2 D} to having a polynomial part of any claimed degree $d>0$ instead of $3(3k+2)$ in the summation.
Then one can see all congruence identities in \eqref{eq:franel 1}--\eqref{eq:franel 3-2} can be proved uniformly by the polynomial reduction method.

Let $F(k)=(-1)^kf_k=(-1)^k\sum_{i=0}^{k}\binom{k}{i}^3$.
By Zeilberger's algorithm, we find that
\begin{equation}\label{eq:L Franel}
L=(k+2)^2\sigma^2+(7k^2+21k+16)\sigma-8(k+1)^2\in \ann F(k).
\end{equation}

\begin{theo}\label{th:Franel cong}
Let $L$ be as in \eqref{eq:L Franel} and $n$ a positive integer. Then
\begin{equation}\label{eq: Franel D}
\sum_{k=0}^{n-1} L^{\ast}(p(k))(-1)^kf_k
=-n^2(p(n-2)F(n)+8p(n-1)F(n-1)).
\end{equation}
for any polynomial $p(k)\in\bZ[k]$.
Here $L^{\ast}$ is the adjoint of $L$.
\end{theo}
\pf
By Equality \eqref{eq:rec ann} and the fact $u_0(0)F(0)+u_1(0)F(1)=0$, we have
\begin{equation}\label{eq:Franel Delta}
\sum_{k=0}^{n-1}L^{\ast}(p(k)) F(k)
=
-\left(u_0(n)F(n)+u_1(n)F(n+1)\right),
\end{equation}
where
$u_0(n)=n^2p(n-2)+(7n^2+7n+2)p(n-1)$ and
$u_1(n)=(n+1)^2p(n-1)$.
As $L\in\ann F(k)$, it is straightforward to check that for any $n\geq 1$ \begin{equation}\label{eq:shift Franel}
(n+1)^2F(n+1)=8n^2F(n-1)-(7n^2+7n+2)F(n).
\end{equation}
Substituting \eqref{eq:shift Franel} into \eqref{eq:Franel Delta} derives \eqref{eq: Franel D}.
\qed

Theorem \ref{th:Franel cong} together with \eqref{eq:franel mod2} lead to the following corollary.
\begin{coro}\label{coro:Franel cong}
Let $L$ be as in \eqref{eq:L Franel}. Then
\begin{equation}\label{eq:Franel cong}
\sum_{k=0}^{n-1} L^{\ast}(p(k))(-1)^kf_k \equiv 0 \pmod {2n^2}
\end{equation}
for any polynomial $p(k)\in\bZ[k]$.
\end{coro}
Since $3(3k+2)=-L^{\ast}(1)$, by \eqref{eq: Franel D} we have
\begin{equation*}
3\sum_{k=0}^{n-1}(3k+2)(-1)^k f_k=
-\sum_{k=0}^{n-1}L^{\ast}(1)(-1)^k f_k=
n^2(F(n)+8F(n-1)),
\end{equation*}
which is exactly \eqref{eq:franel 3k+2 D}.

Equalities \eqref{eq:franel 3-1} and \eqref{eq:franel 3-2} can be proved by \eqref{eq: Franel D} and the observation that
$27k^2(3k+1)=-L^{\ast}(1)-3L^{\ast}(k^2)$ and  $9(9k^3-15k^2-10k)=-4L^{\ast}(1)+9L^{\ast}(k)-3L^{\ast}(k^2).
$
By \eqref{eq:franel 0} and the decompositions
\[
k=-\frac{2}{3}-\frac{1}{9}L^{\ast}(1)
\quad\text{and}\quad
k^2=\frac{10}{27}+\frac{13}{162}L^{\ast}(1)-\frac{1}{18}L^{\ast}(k),
\]
\eqref{eq:franel 1} and \eqref{eq:franel 2} can also be confirmed.

In fact, when $p(k)\in\bZ[k]$ is a polynomial of degree $s\in\bN$, since $L$ in \eqref{eq:L Franel} is not degenerated,
we know $\deg L^{\ast}(p(k))=s+1$.
\begin{coro}
For any positive integer $d$, we can find a polynomial $q(k)\in\bZ[k]$ with $\deg q(k)=d$ such that
\[
\sum_{k=0}^{n-1} q(k)(-1)^kf_k \equiv 0 \pmod {2n^2}.
\]
\end{coro}
The polynomial reduction method also applies to other holonomic sequences.
Let us exhibit with one more example.
The central Delannoy numbers $D_k$ are defined by
\[
D_k=\sum_{i=0}^{k}\binom{k}{i}\binom{k+i}{i}.
\]
Zeilberger's algorithm leads to
\begin{equation}\label{eq:L Delannoy}
L=(k+2)\sigma^2+(-6k-9)\sigma+(k+1)\in \ann D_k.
\end{equation}
Then by a similar argument to the proof of Theorem \ref{th:Franel cong}, we have
\begin{theo}\label{th:Delannoy cong}
Let $L$ be as in \eqref{eq:L Delannoy}. Then
\begin{equation}\label{eq:Delannoy Delta}
\sum_{k=0}^{n-1} L^{\ast}(p(k))D_k
=n(p(n-1)D_{n-1}-p(n-2)D_n).
\end{equation}
for any polynomial $p(k)\in\bZ[k]$.
\end{theo}

When $p(k)=1$, equality \eqref{eq:Delannoy Delta} becomes
\[
\sum_{k=0}^{n-1} (4k+2)D_k
=n(D_{n}-D_{n-1}),
\]
which was first observed by C. Wang (private communication).

\begin{coro}\label{coro:Delannoy cong}
Let $L$ be as in \eqref{eq:L Delannoy}. Then for any polynomial $p(k)\in\bZ[k]$, we have
\begin{equation}\label{eq:Delannoy cong}
\sum_{k=0}^{n-1} L^{\ast}(p(k))D_k \equiv 0 \pmod {n}.
\end{equation}
\end{coro}

\noindent \textbf{Acknowledgments.}
This work was supported by the National Natural Science Foundation of China (No. 12101449, 11701419, 11871067).


\begin{thebibliography}{10}

\bibitem{Abramov1975}
S.A. Abramov.
\newblock The rational component of the solution of a first order linear recurrence relation with rational right hand side.
%{\em {\v{Z}}.Vy{\v{c}}isl. Mat. i Mat. Fiz.}, 15(1975), 1035--1039, 1090.
{\em U.S.S.R. Comput. Math. Math. Phys.}, 15(1975), 216--221.

\bibitem{Abramov1989}
S.A. Abramov.
\newblock Rational solutions of linear differential and difference equations with polynomial coefficients.
{\em U.S.S.R. Comput. Math. Math. Phys.}, 29(1989), 7--12.

\bibitem{Abramov1995}
S.A. Abramov.
\newblock Rational solutions of linear difference and $q$-difference equations with polynomial coefficients.
\newblock In {\em ISSAC '95}, pages 285--289, 1995. ACM.

\bibitem{AbramovvanHoeij}
S.A. Abramov and M. van Hoeij.
\newblock Integration of solutions of linear functional equations.
{\em Integral Transforms Spec. Funct.}, 8(1999), 3--12.

\bibitem{AP2001}
S.A. Abramov and M. Petkov{\v{s}}ek.
\newblock Minimal decomposition of indefinite hypergeometric sums.
\newblock In {\em ISSAC '01}, pages 7--14, 2001. ACM.

\bibitem{AP2002}
S.A. Abramov and M. Petkov{\v{s}}ek.
\newblock Rational normal forms and minimal decompositions of hypergeometric terms.
\newblock {\em J. Symbolic Compt.}, 33(2002), 521--543.

\bibitem{BCCL2010}
A. Bostan, S. Chen, F. Chyzak and Z. Li.
\newblock Complexity of creative telescoping for bivariate rational functions.
\newblock In {\em ISSAC '10}, pages 203--210, 2010. ACM.

%\bibitem{BCCLX2013}
%A. Bostan, S. Chen, F. Chyzak, Z. Li, and G. Xin.
%\newblock Hermite reduction and creative telescoping for hyperexponential functions.
%\newblock In {\em ISSAC '13}, pages 77--84, 2013.

\bibitem{BCLS2018}
A. Bostan, F. Chyzak, P. Lairez and B. Salvy.
\newblock Generalized Hermite reduction, creative telescoping and definite integration of D-finite functions.
\newblock In {\em ISSAC '18}, pages 95--102, 2018. ACM.

\bibitem{BLS2013}
A. Bostan, P. Lairez and B. Salvy.
\newblock Creative telescoping for rational functions using the Griffiths--Dwork method.
\newblock In {\em ISSAC '13},  pages 93--100, 2013. ACM.


\bibitem{CCL2004}
H.H. Chan, S.H. Chan and Z. Liu.
\newblock Domb's numbers and Ramanujan--Sato type series for $1/\pi$.
\newblock {\em Adv. Math.}, 186(2004), 396--410.

\bibitem{CTYZ2011}
H.H. Chan, Y. Tanigawa, Y. Yang and W. Zudilin.
\newblock New analogues of Clausen's identities arising from the theory of modular forms.
\newblock {\em Adv. Math.}, 228(2011), 1294--1314.

\bibitem{CCFFL2015}
S. Chen, F. Chyzak, R. Feng, G. Fu and Z. Li.
\newblock On the existence of telescopers for mixed hypergeometric terms. \newblock {\em J. Symbolic Comput.}, 68(2015), 1--26.

\bibitem{CDZ2019}
S. Chen, L. Du and C. Zhu.
\newblock Existence problem of telescopers for rational functions in three variables: the mixed cases.
\newblock In {\em ISSAC '19},  pages 82--89, 2019. ACM.

\bibitem{CHKK2018}
S. Chen,  M. van Hoeij and M. Kauers.
\newblock Reduction-based creative telescoping for fuchsian D-finite functions.
\newblock {\em J. Symbolic Comput.}, 85(2018), 108--127.

\bibitem{CHHLW}
S. Chen, Q.-H. Hou, H. Huang, G. Labahn and R.-H. Wang.
\newblock Constructing minimal telescopers for rational functions in three discrete variables.
\newblock {\em Adv. in Appl. Math.}, to appear.

\bibitem{CHLW2016}
S. Chen, Q.-H. Hou, G. Labahn and R.-H. Wang.
\newblock Existence problem of telescopers: beyond the bivariate case.
\newblock In {\em ISSAC '16}, pages 167--174, 2016. ACM.

\bibitem{CHKL2015}
S. Chen, H. Huang, M. Kauers and Z. Li.
\newblock A modified Abramov-Petkov\v{s}ek reduction and creative telescoping for hypergeometric terms.
\newblock In {\em ISSAC '15}, pages 117--124, 2015. ACM.

\bibitem{CKK2016}
S. Chen, M. Kauers and C. Koutschan.
\newblock Reduction-based creative telescoping for algebraic functions.
\newblock In {\em ISSAC '16}, pages 175--182, 2016. ACM.

\bibitem{Cooper2012}
S. Cooper.
\newblock Level $10$ analogues of Ramanujan's series for $1/\pi$.
\newblock {\em J. Ramanujan Math. Soc.}, 27(2012), 59--76.

\bibitem{Guo2013}
V.J.W. Guo.
\newblock Proof of two conjectures of Sun on congruences for Franel numbers.
\newblock {\em Integral Transforms Spec. Funct.}, 24(2013), 532--539.

\bibitem{Hermite1872}
C. Hermite.
\newblock Sur l’intégration des fractions rationnelles.
\newblock {\em Ann. Sci. École Norm. Sup. (2)}, 1(1872), 215--218.

\bibitem{Hoeven2018}
J. van der Hoeven.
\newblock Creative telescoping using reductions.
\newblock {\em Preprint:hal-01773137v2}, June 2018.

\bibitem{Hoeven2017}
J. van der Hoeven.
\newblock Constructing reductions for creative telescoping.
\newblock {\em  Appl. Algebra Engrg. Comm. Comput.}, 32(2021), 575--602.

\bibitem{HouLi2021}
Q.-H. Hou and G.-J. Li.
\newblock Gosper summability of rational multiples of hypergeometric terms.
\newblock {\em J. Difference Equ. Appl.}, 27(2021), 1723--1733.

\bibitem{HouLiu2021}
Q.-H. Hou and K. Liu.
\newblock Congurences and telescopings of $P$-recursive sequences.
\newblock {\em J. Difference Equ. Appl.}, 27(2021), 686--697.

\bibitem{HouMuZeil2021}
Q.-H. Hou, Y.-P. Mu and D. Zeilberger.
\newblock Polynomial reduction and supercongruences.
\newblock {\em J. Symbolic Comput.}, 103(2021), 127--140.

\bibitem{Huang2016}
H. Huang.
\newblock New bounds for creative telescoping.
\newblock In {\em ISSAC '16}, pages 279–286, 2016. ACM.

\bibitem{KauersPaule2011}
M. Kauers and P. Paule.
\newblock {\em The Concrete Tetrahedron}.
\newblock Springer Wien, 2011.

\bibitem{Lairez2016}
P. Lairez.
\newblock Computing periods of rational integrals.
\newblock {\em Math. Comp.}, 85(2016), 1719--1752.

\bibitem{O1845}
M. Ostrogradsky.
\newblock De l'integration des fractions rationelles.
\newblock {\em Bull. de la Classe Physico--Mathématique de l’Acad. Impériale des Sciences de St.-Pétersbourg}, 4(1845), 145--167, 286--300.

\bibitem{Rogers2009}
M.D. Rogers.
\newblock New $_5F_4$ hypergeometric transformations, three-variable Mahler measures, and formulas for $1/\pi$.
\newblock {\em Ramanujan J.}, 18(2009), 327--340.

\bibitem{Sun2013}
Z.-W. Sun.
\newblock Congruences for Franel numbers.
\newblock {\em Adv. in Appl. Math.}, 51(2013), 524--535.

\bibitem{Sun2021}
Z.-W. Sun.
\newblock New type series for powers of $\pi$.
\newblock arXiv:2110.03651.

\bibitem{Trager1984}
B.M. Trager.
\newblock Integration of Algebraic Functions.
\newblock {\em Ph.D. thesis, MIT}, 1984.

\bibitem{WangSun2019}
C. Wang and Z.-W. Sun.
\newblock Divisility results on Franel numbers and related polynomials.
\newblock {\em Int. J. Number Theory}, 15(2019), 433--444.

\bibitem{WZ2022}
R.-H. Wang and M.X.X. Zhong.
\newblock $q$-Rational reduction and $q$-analogues of series for $\pi$.
\newblock arXiv:2203.16047.

\bibitem{WZ1990}
H.S. Wilf and D. Zeilberger.
\newblock An algorithmic proof theory for hypergeometric (ordinary and ``$q$'') multisum/integral identities.
\newblock {\em Invent. Math.}, 108(1992), 575--633.

\bibitem{Zeilberger1990c}
D. Zeilberger.
\newblock A fast algorithm for proving terminating hypergeometric identities.
\newblock {\em Discrete Math.}, 80(1990), 207--211.

\bibitem{Zeil1990}
D. Zeilberger.
\newblock A holonomic systems approach to special function identities.
\newblock {\em J. Comput. Appl. Math.}, 32(1990), 321–368.

\bibitem{Zeilberger1991}
D. Zeilberger.
\newblock The method of creative telescoping.
\newblock {\em J. Symbolic Comput.}, 11(1991), 195--204.

\end{thebibliography}
\end{document}